\documentclass{article}
\usepackage{amsfonts}
\usepackage{amsmath}

\input{tcilatex}
\begin{document}

\begin{center}
\bigskip

\textbf{CLASSIFICATION OF HOMOTHETICAL HYPERSURFACES AND ITS APPLICATIONS TO
PRODUCTION FUNCTIONS IN ECONOMICS}

\bigskip

\textbf{Muhittin Evren AYDIN, Mahmut ERGUT}

Department of Mathematics, Faculty of Science, Firat University

23200 Elazig, Turkey, meaydin@firat.edu.tr, mergut@firat.edu.tr
\end{center}

\textbf{Abstract. }In this paper, we completely classify the homothetical
hypersurfaces having null Gauss-Kronocker curvature in a Euclidean $\left(
n+1\right) -$space $%
\mathbb{R}
^{n+1}.$ Several applications to the production functions in economics are
also given.

\bigskip


\bigskip


\begin{center}
\textbf{1. Introduction and preliminaries on hypersurfaces}
\end{center}

A hypersurface $M^{n}$ of a Euclidean $\left( n+1\right) -$space $%
\mathbb{R}
^{n+1}$ is called a \textit{homothetical hypersurface} if it is the graph of
a function of the form: 
\begin{equation}
f\left( x_{1},...,x_{n}\right) =f_{1}\left( x_{1}\right) \times ...\times
f_{n}\left( x_{n}\right) ,  \tag{1.1}
\end{equation}%
where $f_{1},...,f_{n}$ are functions of class $C^{\infty }$ \cite{25}. We
call $f_{1},...,f_{n}$ the \textit{components} of $f,$ and denote the
homothetical hypersurface $M^{n}$ by a pair $\left( M^{n},f\right) $.

Homothetical hypersurfaces have been studied by many authors based on
minimality property of these hypersurfaces \cite{18,20,24,25,27}.

G.E. V\^{\i}lcu and A.D. V\^{\i}lcu \cite{28,29} established an interesting
link between some fundamental notions in the theory of production functions
and the differential geometry of hypersurfaces. For further study of
production hypersurfaces, we refer the reader to B.-Y. Chen's series of
interesting papers on homogeneous production functions, quasi-sum production
models, and homothetic production functions [5, 6, 8-15] and X. Wang and Y.
Fu's paper \cite{30}.

Let $M^{n}$ be a hypersurface of a Euclidean $\left( n+1\right) -$space $%
\mathbb{%
\mathbb{R}
}^{n+1}$. For general references on the geometry of hypersurfaces see \cite%
{7,17,19}.

The \textit{Gauss map} $\nu :M^{n}\longrightarrow S^{n+1}$ maps $M^{n}$ to
the unit hypersphere $S^{n}$ of $\mathbb{%
\mathbb{R}
}^{n+1}.$ The differential $d\nu $ of the Gauss map $\nu $ is known as 
\textit{shape operator} or Weingarten map. Denote by $T_{p}M^{n}$ the
tangent space of $M^{n}$ at the point $p\in M^{n}.$ Then, for $v,w\in
T_{p}M^{n},$ the shape operator $A_{p}$ at the point \ $p\in M^{n}$ is
defined by%
\begin{equation*}
g\left( A_{p}\left( v\right) ,w\right) =g\left( d\nu \left( v\right)
,w\right) ,
\end{equation*}%
where $g$ is the induced metric tensor on $M^{n}$ from the Euclidean metric
on $\mathbb{%
\mathbb{R}
}^{n+1}.$

The determinant of the shape operator $A_{p}$ is called the \textit{%
Gauss-Kronocker curvature}. A hypersurface having null Gauss-Kronecker
curvature is said to be \textit{developable}. In this case the hypersurface
can be flattened onto a hyperplane without distortion. We remark that
cylinders and cones are examples of developable surfaces, but the spheres
are not under any metric.

For a given function $f=f\left( x_{1},...,x_{n}\right) ,$ the graph of $f$
is the non-parametric hypersurface of $\mathbb{%
\mathbb{R}
}^{n+1}$ defined by%
\begin{equation}
\varphi \left( \mathbf{x}\right) =\left( x_{1},...,x_{n},f\left( \mathbf{x}%
\right) \right)  \tag{1.2}
\end{equation}%
for $\mathbf{x=}\left( x_{1},...,x_{n}\right) \in 
\mathbb{R}
^{n}$.

Let us put%
\begin{equation}
\omega =\sqrt{1+\sum_{i=1}^{n}\left( \frac{\partial f}{\partial x_{i}}%
\right) ^{2}}.  \tag{1.3}
\end{equation}%
The Gauss-Kronecker curvature of the graph of $f$ is%
\begin{equation}
G=\frac{\det \left( \mathcal{H}\left( f\right) \right) }{\omega ^{n+2}}, 
\tag{1.4}
\end{equation}%
where $\mathcal{H}\left( f\right) $ is the Hessian matrix of $f,$ that is,
the square matrix $\left( f_{x_{i}x_{j}}\right) $ of second-order partial
derivatives of $f.$

In this paper, we completely classify homothetical hypersurfaces having null
Gauss-Kronocker curvature. Several applications to production models in
economics are also given.

\begin{center}
\textbf{2. Production models in economics}
\end{center}

In economics, a \textit{production function} is a mathematical expression
which denotes the physical relations between the output generated of a firm,
an industry or an economy and inputs that have been used. Explicitly, a
production function is a map which has non-vanishing first derivatives
defined by%
\begin{equation*}
f:%
\mathbb{R}
_{+}^{n}\longrightarrow 
\mathbb{R}
_{+},\text{ \ }f=f\left( x_{1},x_{2},...,x_{n}\right) ,
\end{equation*}%
where $f$ is the quantity of output, $n$ are the number of inputs and $%
x_{1},x_{2},...,x_{n}$ are the inputs.

A production function $f\left( x_{1},x_{2},...,x_{n}\right) $ is said to be 
\textit{homogeneous of degree }$p$\textit{\ }or\textit{\ }$p-$\textit{%
homogenous} if%
\begin{equation}
f\left( tx_{1},...,tx_{n}\right) =t^{p}f\left( x_{1},...,x_{n}\right) 
\tag{2.1}
\end{equation}%
holds for each $t\in 
\mathbb{R}
_{+}$ for which $\left( 2.1\right) $ is defined. A homogeneous function of
degree one is called \textit{linearly homogeneous.} If $h>1$, the function
exhibits increasing return to scale, and it exhibits decreasing return to
scale if $h<1$. If it is homogeneous of degree 1, it exhibits constant
return to scale \cite{8}.

In 1928, C. W. Cobb and P. H. Douglas introduces \cite{16} a famous
two-factor production function%
\begin{equation*}
Y=bL^{k}C^{1-k},
\end{equation*}%
where $b$ presents the total factor productivity, $Y$ the total production, $%
L$ the labor input and $C$ the capital input. This function is nowadays
called \textit{Cobb-Douglas production function}. In its generalized form
the Cobb-Douglas production function\textit{\ }may be expressed as%
\begin{equation*}
f\left( \mathbf{x}\right) =\gamma x_{1}^{\alpha _{1}}...x_{n}^{\alpha _{n}},
\end{equation*}%
where $\gamma $ is a positive constant and $\alpha _{1},...,\alpha _{n}$ are
nonzero constants.

In 1961, K. J. Arrow, H. B. Chenery, B. S. Minhas and R. M. Solow \cite{2}
introduced a two-factor production function given by%
\begin{equation*}
Q=F\cdot \left( aK^{r}+\left( 1-a\right) L^{r}\right) ^{\frac{1}{r}},
\end{equation*}%
where $Q$ is the output, $F$ the factor productivity, $a$ the share
parameter, $K$ and $L$ the primary production factors, $r=\left( s-1\right)
/s,$ and $s=1/\left( 1-r\right) $ is the elasticity of substitution.

The \textit{generalized ACMS production function of }$n$ variables is given
by%
\begin{equation*}
f\left( x_{1},...,x_{n}\right) =\gamma \left( \beta
_{1}^{p}x_{1}^{p}+...+\beta _{n}^{p}x_{n}^{p}\right) ^{\frac{d}{p}},
\end{equation*}%
where $\rho \neq 0,$ $\rho <1,$ $d,\gamma >0$ and $\beta _{i}>0$ for all $%
i=1,...,n$.

\textit{A homothetic function} is a production function of the form:%
\begin{equation*}
f\left( x_{1},...,x_{n}\right) =F\left( h\left( x_{1},...,x_{n}\right)
\right) ,
\end{equation*}%
where $h\left( x_{1},...x_{n}\right) $ is homogeneous function of arbitrary
given degree and $F$ is a monotonically increasing function \cite{12,15,21}.

A homothetic production function of form%
\begin{equation*}
f\left( \mathbf{x}\right) =F\left( \sum_{i=1}^{n}\beta _{i}^{\rho
}x_{i}^{\rho }\right) ^{\frac{d}{\rho }}\left( \text{resp., }f\left( \mathbf{%
x}\right) =F\left( x_{1}^{\alpha _{1}}...x_{n}^{\alpha _{n}}\right) \right)
\end{equation*}%
is called a\textit{\ homothetic generalized ACMS production function}
(resp., a \textit{homothetic generalized} \textit{Cobb-Douglas production
function}) \cite{11}.

The most common quantitative indices of production factor substitutability
are forms of the elasticity of substitution. R.G.D. Allen and J.R. Hicks 
\cite{1} suggested two generalizations of Hicks' original two variable
elasticity concept.

The first concept, called Hicks elasticity of substitution\textit{,} is
defined as follows.

Let $f\left( x_{1},...,x_{n}\right) $ be a production function. Then \textit{%
Hicks elasticity of substitution} of the $i-$th production variable with
respect to the $j-$th production variable is given by%
\begin{equation}
H_{ij}\left( \mathbf{x}\right) =-\frac{\dfrac{1}{x_{i}f_{x_{i}}}+\dfrac{1}{%
x_{j}f_{x_{j}}}}{\dfrac{f_{x_{i}x_{i}}}{\left( f_{x_{i}}\right) ^{2}}-\dfrac{%
2f_{x_{i}x_{j}}}{f_{x_{i}}f_{x_{j}}}+\dfrac{f_{x_{j}x_{j}}}{\left(
f_{x_{j}}\right) ^{2}}}\text{ \ }\left( \mathbf{x\in }%
\mathbb{R}
_{+}^{n},\text{ }i,j=1,...,n,\text{ }i\neq j\right) ,  \tag{2.2}
\end{equation}%
where $f_{x_{i}}=\partial f/\partial x_{i},$ $f_{x_{i}x_{j}}=\partial
^{2}f/\partial x_{i}\partial x_{j}.$

A production function $f$ is said to satisfy the constant Hicks elasticity
of substitution property if there is a nonzero constant $\sigma \in 
\mathbb{R}
$ such that%
\begin{equation}
H_{ij}\left( \mathbf{x}\right) =\sigma ,\text{ for }\mathbf{x\in }%
\mathbb{R}
_{+}^{n}\text{ and }1\leq i\neq j\leq n.  \tag{2.3}
\end{equation}

L. Losonczi \cite{22} classified homogeneous production functions of $2$
variables, having constant Hicks elasticiy of substitution. Then, the
classification of L. Losonczi was extended to $n$ variables by B-Y. Chen 
\cite{13}.

The second concept, investigated by R.G.D. Allen and H. Uzawa \cite{26}$,$
is the following:

Let $f$ be a production function. Then \textit{Allen elasticity of
substitution} of the $i-$th production variable with respect to the $j-$th
production variable is defined by

\begin{equation}
A_{ij}\left( \mathbf{x}\right) =-\frac{%
x_{1}f_{x_{1}}+x_{2}f_{x_{2}}+...+x_{n}f_{x_{n}}}{x_{i}x_{j}}\frac{D_{ij}}{%
\det \left( \mathcal{H}^{B}\left( f\right) \right) }\text{ \ }\left( \mathbf{%
x}=\left( x_{1},...,x_{n}\right) \mathbf{\in }%
\mathbb{R}
_{+}^{n},\text{ }i,j=1,...,n,\text{ }i\neq j\right) ,  \tag{2.4}
\end{equation}%
where $D$ is the determinant of the bordered Hessian matrix%
\begin{equation}
\mathcal{H}^{B}\left( f\right) =%
\begin{pmatrix}
0 & f_{x_{1}} & ... & f_{x_{n}} \\ 
f_{x_{1}} & f_{x_{1}x_{1}} & ... & f_{x_{1}x_{n}} \\ 
\vdots & \vdots & ... & \vdots \\ 
f_{x_{n}} & f_{x_{n}x_{1}} & ... & f_{x_{n}x_{n}}%
\end{pmatrix}%
,  \tag{2.5}
\end{equation}%
and $D_{ij}$ is the co-factor of the element $f_{x_{i}x_{j}}$ in the
determinant $\det \left( \mathcal{H}^{B}\left( f\right) \right) $ ($\det
\left( \mathcal{H}^{B}\left( f\right) \right) \neq 0$ is assumed). The
authors call the bordered Hessian matrix $\mathcal{H}^{B}\left( f\right) $
by \textit{Allen's matrix} and $\det \mathcal{H}^{B}\left( f\right) $ by 
\textit{Allen determinant} in \cite{3,4}.

It is a simple calculation to show that in case of two variables Hicks
elasticity of substitution coincides with Allen elasticity of substitution.

\begin{center}
\textbf{3. Classification of homothetical hypersurfaces}
\end{center}

Throughout this article, we assume that the functions $f_{1}\left(
x_{1}\right) ,...,f_{n}\left( x_{n}\right) $ are real valued functions and
have non-vanishing first derivatives.

The following provides an useful formula for Hessian determinant of a
function of the form $\left( 1.1\right) $

\bigskip

\textbf{Lemma 3.1. }\textit{The determinant of the Hessian matrix of the
function }$f\left( \mathbf{x}\right) =f_{1}\left( x_{1}\right) \times
...\times f_{n}\left( x_{n}\right) $\textit{\ is given by}%
\begin{equation*}
\det \left( H\left( f\right) \right) =\left( f\right) ^{n}\left[ \frac{%
f_{1}^{\prime \prime }}{f_{1}}\dprod\limits_{i=2}^{n}\left( \frac{%
f_{i}^{\prime }}{f_{i}}\right) ^{\prime }+\left( \frac{f_{1}^{\prime }}{f_{1}%
}\right) ^{\prime }\sum_{i=2}^{n}\left( \frac{f_{2}^{\prime }}{f_{2}}\right)
^{\prime }...\left( \frac{f_{i-1}^{\prime }}{f_{i-1}}\right) ^{\prime
}\left( \frac{f_{i}^{\prime }}{f_{i}}\right) ^{2}\left( \frac{%
f_{i+1}^{\prime }}{f_{i+1}}\right) ^{\prime }...\left( \frac{f_{n}^{\prime }%
}{f_{n}}\right) ^{\prime }\right] ,
\end{equation*}%
where $f_{i}^{\prime }=\frac{df}{dx_{i}},$ $f_{i}^{\prime \prime }=\frac{%
d^{2}f}{dx_{i}^{2}}$ for all $i\in \left\{ 1,...,n\right\} .$

\bigskip

\textbf{Proof. }Let $f$ be a twice differentiable function given by%
\begin{equation}
f\left( \mathbf{x}\right) =f_{1}\left( x_{1}\right) \times ...\times
f_{n}\left( x_{n}\right)  \tag{3.1}
\end{equation}%
for $\mathbf{x}=\left( x_{1},...,x_{n}\right) \mathbf{\in 
\mathbb{R}
}^{n}.$ It follows from $\left( 3.1\right) $ that%
\begin{equation}
f_{x_{i}}=\frac{f_{i}^{\prime }}{f_{i}}f,\text{ }f_{x_{i}x_{j}}=\frac{%
f_{i}^{\prime }f_{j}^{\prime }}{f_{i}f_{j}}f,\text{ }f_{x_{i}x_{i}}=\frac{%
f_{i}^{\prime \prime }}{f_{i}}f,\text{ }1\leq i\neq j\leq n.  \tag{3.2}
\end{equation}%
By using $\left( 3.2\right) ,$ the Hessian determinant of the function $f$
is 
\begin{equation}
\det \left( H\left( f\right) \right) =\left( f\right) ^{n}%
\begin{vmatrix}
\dfrac{f_{1}^{\prime \prime }}{f_{1}} & \dfrac{f_{1}^{\prime }f_{2}^{\prime }%
}{f_{1}f_{2}} & \dfrac{f_{1}^{\prime }f_{3}^{\prime }}{f_{1}f_{3}} & ... & 
\dfrac{f_{1}^{\prime }f_{n}^{\prime }}{f_{1}f_{n}} \\ 
\dfrac{f_{1}^{\prime }f_{2}^{\prime }}{f_{1}f_{2}} & \dfrac{f_{2}^{\prime
\prime }}{f_{2}} & \dfrac{f_{2}^{\prime }f_{3}^{\prime }}{f_{2}f_{3}} & ...
& \dfrac{f_{2}^{\prime }f_{n}^{\prime }}{f_{2}f_{n}} \\ 
\dfrac{f_{1}^{\prime }f_{3}^{\prime }}{f_{1}f_{3}} & \dfrac{f_{2}^{\prime
}f_{3}^{\prime }}{f_{2}f_{3}} & \dfrac{f_{3}^{\prime \prime }}{f_{3}} & ...
& \dfrac{f_{3}^{\prime }f_{n}^{\prime }}{f_{3}f_{n}} \\ 
\vdots & \vdots & \vdots & ... & \vdots \\ 
\dfrac{f_{1}^{\prime }f_{n}^{\prime }}{f_{1}f_{n}} & \dfrac{f_{2}^{\prime
}f_{n}^{\prime }}{f_{2}f_{n}} & \dfrac{f_{3}^{\prime }f_{n}^{\prime }}{%
f_{3}f_{n}} & ... & \dfrac{f_{n}^{\prime \prime }}{f_{n}}%
\end{vmatrix}%
.  \tag{3.3}
\end{equation}%
Now we apply Gauss elimination method for the determinant from the formula $%
\left( 3.3\right) $. We replace the second column by second column minus $%
\left( \frac{f_{1}f_{2}^{\prime }}{f_{1}^{\prime }f_{2}}\right) $ times
first column; then we derive%
\begin{equation*}
\det \left( H\left( f\right) \right) =\left( f\right) ^{n}%
\begin{vmatrix}
\dfrac{f_{1}^{\prime \prime }}{f_{1}} & -\dfrac{f_{1}f_{2}^{\prime }}{%
f_{1}^{\prime }f_{2}}\left( \dfrac{f_{1}^{\prime }}{f_{1}}\right) ^{\prime }
& \dfrac{f_{1}^{\prime }f_{3}^{\prime }}{f_{1}f_{3}} & ... & \dfrac{%
f_{1}^{\prime }f_{n}^{\prime }}{f_{1}f_{n}} \\ 
\dfrac{f_{1}^{\prime }f_{2}^{\prime }}{f_{1}f_{2}} & \left( \dfrac{%
f_{2}^{\prime }}{f_{2}}\right) ^{\prime } & \dfrac{f_{2}^{\prime
}f_{3}^{\prime }}{f_{2}f_{3}} & ... & \dfrac{f_{2}^{\prime }f_{n}^{\prime }}{%
f_{2}f_{n}} \\ 
\dfrac{f_{1}^{\prime }f_{3}^{\prime }}{f_{1}f_{3}} & 0 & \dfrac{%
f_{3}^{\prime \prime }}{f_{3}} & ... & \frac{f_{3}^{\prime }f_{n}^{\prime }}{%
f_{2}f_{n}} \\ 
\vdots & \vdots & \vdots & ... & \vdots \\ 
\dfrac{f_{1}^{\prime }f_{n}^{\prime }}{f_{1}f_{n}} & 0 & \dfrac{%
f_{3}^{\prime }f_{n}^{\prime }}{f_{3}f_{n}} & ... & \dfrac{f_{n}^{\prime
\prime }}{f_{n}}%
\end{vmatrix}%
\end{equation*}%
By similar elementary transformations, we get%
\begin{equation}
\det \left( H\left( f\right) \right) =\left( f\right) ^{n}%
\begin{vmatrix}
\dfrac{f_{1}^{\prime \prime }}{f_{1}} & -\dfrac{f_{1}f_{2}^{\prime }}{%
f_{1}^{\prime }f_{2}}\left( \dfrac{f_{1}^{\prime }}{f_{1}}\right) ^{\prime }
& -\dfrac{f_{1}f_{3}^{\prime }}{f_{1}^{\prime }f_{3}}\left( \dfrac{%
f_{1}^{\prime }}{f_{1}}\right) ^{\prime } & ... & -\dfrac{f_{1}f_{n}^{\prime
}}{f_{1}^{\prime }f_{n}}\left( \dfrac{f_{1}^{\prime }}{f_{1}}\right)
^{\prime } \\ 
\dfrac{f_{1}^{\prime }f_{2}^{\prime }}{f_{1}f_{2}} & \left( \dfrac{%
f_{2}^{\prime }}{f_{2}}\right) ^{\prime } & 0 & ... & 0 \\ 
\dfrac{f_{1}^{\prime }f_{3}^{\prime }}{f_{1}f_{3}} & 0 & \left( \dfrac{%
f_{3}^{\prime }}{f_{3}}\right) ^{\prime } & ... & 0 \\ 
\vdots & \vdots & \vdots & ... & \vdots \\ 
\dfrac{f_{1}^{\prime }f_{n}^{\prime }}{f_{1}f_{n}} & 0 & 0 & ... & \dfrac{%
f_{n}^{\prime \prime }}{f_{n}}%
\end{vmatrix}%
.  \tag{3.4}
\end{equation}%
After calculating the determinant from the formula $\left( 3.4\right) ,$ we
finally obtain%
\begin{equation*}
\det \left( H\left( f\right) \right) =\left( f\right) ^{n}\left[ \frac{%
f_{1}^{\prime \prime }}{f_{1}}\left( \frac{f_{2}^{\prime }}{f_{2}}\right)
^{\prime }...\left( \frac{f_{n}^{\prime }}{f_{n}}\right) ^{\prime }+\left( 
\frac{f_{1}^{\prime }}{f_{1}}\right) ^{\prime }\sum_{i=2}^{n}\left( \frac{%
f_{2}^{\prime }}{f_{2}}\right) ^{\prime }...\left( \frac{f_{i-1}^{\prime }}{%
f_{i-1}}\right) ^{\prime }\left( \frac{f_{i}^{\prime }}{f_{i}}\right)
^{2}\left( \frac{f_{i+1}^{\prime }}{f_{i+1}}\right) ^{\prime }...\left( 
\frac{f_{n}^{\prime }}{f_{n}}\right) ^{\prime }\right] .
\end{equation*}

Next result completely classifies the homothetical hypersurfaces having null
Gauss-Kronocker curvature.

\bigskip

\textbf{Theorem\textbf{\ }3.1. }\textit{Let }$\left( M^{n},f\right) $\textit{%
\ be a homothetical hypersurface in }$%
\mathbb{R}
^{n+1}.$\textit{\ }$\left( M^{n},f\right) $\textit{\ has null
Gauss-Kronocker curvature if and only if it is parametrized by one of the
following}

\textit{(a) }$\varphi \left( \mathbf{x}\right) =\left(
x_{1},...,x_{n},f_{1}\left( x_{1}\right) \times \gamma e^{\lambda
_{2}x_{2}+\lambda _{3}x_{3}}\times ...\times f_{n}\left( x_{n}\right)
\right) $\textit{\ for nonzero constants }$\gamma ,\lambda _{2},\lambda
_{3}; $

\textit{(b) }$\varphi \left( \mathbf{x}\right) =\left(
x_{1},...,x_{n},\gamma \left( x_{1}+\beta _{1}\right) ^{\alpha _{1}}\times
...\times \left( x_{n}+\beta _{n}\right) ^{\alpha _{n}}\right) $\textit{,
where }$\beta _{1},...,\beta _{n}$\textit{\ are some constants }$\gamma
,\alpha _{1},...,\alpha _{n}$\textit{\ nonzero constants such that }$%
\sum_{i=1}^{n}\alpha _{i}=1.$

\bigskip

\textbf{Proof. }Let $\left( M^{n},f\right) $\textit{\ be a homothetical
hypersurface in }$%
\mathbb{R}
^{n+1}$ parametrized by%
\begin{equation*}
\varphi \left( \mathbf{x}\right) =\left( x_{1},...,x_{n},f_{1}\left(
x_{1}\right) \times ...\times f_{n}\left( x_{n}\right) \right) .
\end{equation*}%
Assume that $\left( M^{n},f\right) $ has null Gauss-Kronocker curvature. It
follows from $\left( 1\text{.}4\right) $ that $\det \left( H\left( f\right)
\right) =0.$ Hence by Lemma 3.1, we get%
\begin{equation}
\frac{f_{1}^{\prime \prime }}{f_{1}}\dprod\limits_{i=2}^{n}\left( \frac{%
f_{i}^{\prime }}{f_{i}}\right) ^{\prime }+\left( \frac{f_{1}^{\prime }}{f_{1}%
}\right) ^{\prime }\sum_{i=2}^{n}\left( \frac{f_{2}^{\prime }}{f_{2}}\right)
^{\prime }...\left( \frac{f_{i-1}^{\prime }}{f_{i-1}}\right) ^{\prime
}\left( \frac{f_{i}^{\prime }}{f_{i}}\right) ^{2}\left( \frac{%
f_{i+1}^{\prime }}{f_{i+1}}\right) ^{\prime }...\left( \frac{f_{n}^{\prime }%
}{f_{n}}\right) ^{\prime }=0.  \tag{3.5}
\end{equation}%
For the equation $\left( 3.5\right) $ we have two cases:

\bigskip

\textbf{Case (i): }At least one of $\left( \frac{f_{1}^{\prime }}{f_{1}}%
\right) ^{\prime },...,\left( \frac{f_{n}^{\prime }}{f_{n}}\right) ^{\prime
} $ vanishes. Without loss of generality, we may assume that 
\begin{equation}
\left( \frac{f_{2}^{\prime }}{f_{2}}\right) ^{\prime }=0.  \tag{3.6}
\end{equation}%
Thus from $\left( 3.5\right) $ we have%
\begin{equation}
\left( \frac{f_{1}^{\prime }}{f_{1}}\right) ^{\prime }\left( \frac{%
f_{2}^{\prime }}{f_{2}}\right) ^{2}\left( \frac{f_{3}^{\prime }}{f_{3}}%
\right) ^{\prime }...\left( \frac{f_{n}^{\prime }}{f_{n}}\right) ^{\prime
}=0.  \tag{3.7}
\end{equation}%
Without loss of generality, we may assume from $\left( 3.7\right) $ that%
\begin{equation}
\left( \frac{f_{3}^{\prime }}{f_{3}}\right) ^{\prime }=0.  \tag{3.8}
\end{equation}%
By solving $\left( 3.6\right) $ and $\left( 3.8\right) ,$ we conclude the
following%
\begin{equation*}
f_{2}\left( x_{2}\right) =\gamma _{2}e^{\lambda _{2}x_{2}},\text{ \ \ }%
f_{3}\left( x_{3}\right) =\gamma _{3}e^{\lambda _{3}x_{3}}
\end{equation*}%
for nonzero constants $\gamma _{2},\gamma _{3},\lambda _{2},\lambda _{3}.$
This gives the statement (a) of the theorem.

\bigskip

\textbf{Case (ii): }$\left( \frac{f_{1}^{\prime }}{f_{1}}\right) ^{\prime
},...,\left( \frac{f_{n}^{\prime }}{f_{n}}\right) ^{\prime }$ are nonzero.
Then from $\left( 3.5\right) ,$ by dividing with the product $\left( \frac{%
f_{1}^{\prime }}{f_{1}}\right) ^{\prime }\times ...\times \left( \frac{%
f_{n}^{\prime }}{f_{n}}\right) ^{\prime }$, we write%
\begin{equation}
\frac{\frac{f_{1}^{\prime \prime }}{f_{1}}}{\left( \frac{f_{1}^{\prime }}{%
f_{1}}\right) ^{\prime }}+\left\{ \frac{\left( \frac{f_{2}^{\prime }}{f_{2}}%
\right) ^{2}}{\left( \frac{f_{2}^{\prime }}{f_{2}}\right) ^{\prime }}+...+%
\frac{\left( \frac{f_{n}^{\prime }}{f_{n}}\right) ^{2}}{\left( \frac{%
f_{n}^{\prime }}{f_{n}}\right) ^{\prime }}\right\} =0.  \tag{3.9}
\end{equation}

We divide the proof of case (ii) into two cases.

\bigskip

\textbf{Case (ii.a): }Taking the partial derivative of $\left( 3.9\right) $
with respect to $x_{i}$ for $i=2,3,...,n,$ we have%
\begin{equation}
\left( \frac{f_{i}^{\prime }}{f_{i}}\right) \left( \frac{f_{2}^{\prime }}{%
f_{2}}\right) ^{\prime \prime }=2\left[ \left( \frac{f_{i}^{\prime }}{f_{i}}%
\right) ^{\prime }\right] ^{2}.  \tag{3.10}
\end{equation}%
By solving $\left( 3.10\right) ,$ we find%
\begin{equation}
f_{i}\left( x_{i}\right) =\gamma _{i}\left( x_{i}+\beta _{i}\right) ^{\alpha
_{i}},\text{ }2\leq i\leq n  \tag{3.11}
\end{equation}%
for some nonzero constants $\alpha _{i},\gamma _{i}$ and some constants $%
\beta _{i}.$

\bigskip

\textbf{Case (ii.b): }Taking the partial derivative of $\left( 3.9\right) $
with respect to $x_{1},$ we get

\begin{equation*}
\frac{f_{1}^{\prime \prime \prime }\left( x_{1}\right) }{f_{1}^{\prime
\prime }\left( x_{1}\right) }+\frac{f_{1}^{\prime }\left( x_{1}\right) }{%
f_{1}\left( x_{1}\right) }=2\frac{f_{1}^{\prime \prime }\left( x_{1}\right) 
}{f_{1}^{\prime }\left( x_{1}\right) },
\end{equation*}%
which implies that%
\begin{equation}
f_{1}\left( x_{1}\right) f_{1}^{\prime \prime }\left( x_{1}\right) =\tau
\left( f_{1}^{\prime }\left( x_{1}\right) \right) ^{2}  \tag{3.12}
\end{equation}%
for some nonzero constant $\tau .$

Now, we divide the proof of case(ii.b) into two cases based on the value of $%
\tau .$

\bigskip

\textbf{Case (ii.b.1): }$\tau =1.$ This case is not possible because of $%
\left( \frac{f_{1}^{\prime }}{f_{1}}\right) ^{\prime },...,\left( \frac{%
f_{n}^{\prime }}{f_{n}}\right) ^{\prime }$ are nonzero.

\bigskip

\textbf{Case (ii.b.2): }$\tau \neq 1.$ After solving $\left( 3.12\right) $,
we derive%
\begin{equation}
f_{1}\left( x_{1}\right) =\gamma _{1}\left( x_{1}+\beta _{1}\right) ^{-\frac{%
1}{\tau -1}}.  \tag{3.13}
\end{equation}%
By substituting $\left( 3.11\right) $ and $\left( 3.13\right) $ into $\left(
3.9\right) ,$ we deduce that $\alpha _{2}+...+\alpha _{n}=\tau /\tau -1.$
Therefore we obtain case (b) of the theorem.

\bigskip

Conversely, it is direct to verify all homotetical hypersurfaces
parametrized by cases (a) and (b) have null Gauss-Kronocker curvature.

\begin{center}
\textbf{4. Applications to Cobb-Douglas production functions}
\end{center}

Geometric representation of the generalized Cobb-Douglas production is given
by the hypersurface $\varphi :%
\mathbb{R}
_{+}^{n}\longrightarrow 
\mathbb{R}
_{+}^{n+1},$%
\begin{equation*}
\varphi \left( \mathbf{x}\right) =\left( x_{1},...,x_{n},\gamma
x_{1}^{\alpha _{1}}...x_{n}^{\alpha _{n}}\right) ,
\end{equation*}%
which is called the \textit{Cobb-Douglas hypersurface} \cite{28}. G. E.Vilcu 
\cite{28} proved that a generalized Cobb-Douglas hypersurface is developable
if and only if it has constant return to scale, i.e, $\sum_{i=1}^{n}\alpha
_{i}=1$. Thus we have the following as a consequence of Theorem 4.1:

\bigskip

\textbf{Corollary 4.1. }\textit{Let }$\left( M^{n},f\right) $\textit{\ be a
homothetical hypersurface in }$%
\mathbb{R}
_{+}^{n+1}$ such that all components of $f$ satisfy $\left( \frac{%
f_{i}^{\prime }}{f_{i}}\right) ^{\prime }\neq 0.$\textit{\ }$\left(
M^{n},f\right) $\textit{\ has null Gauss-Kronocker curvature if and only if,
up to constants, it is a generalized Cobb-Douglas hypersurface having
constant return to scale.}

\bigskip

On the other hand, assume that $h_{i}:%
\mathbb{R}
_{+}\longrightarrow 
\mathbb{R}
$ $\left( i=1,...,n\right) $ and $F:I\subset 
\mathbb{R}
\longrightarrow 
\mathbb{R}
_{+}$ are non-vanishing differentiable functions having nonzero first
derivatives. Then for $h_{1}\left( x_{1}\right) \times ...\times h_{n}\left(
x_{n}\right) \in I,$ we have the following composite function%
\begin{equation}
f\left( x_{1},...,x_{n}\right) =F\left( h_{1}\left( x_{1}\right) \times
...\times h_{n}\left( x_{n}\right) \right) .  \tag{4.1}
\end{equation}

The authors obtained following result in \cite{4}:

\bigskip

\textbf{Theorem 4.1. }\cite{4} \textit{Let }$F\left( u\right) $\textit{\ be
a twice differentiable function with }$F^{\prime }\left( u\right) \neq 0$%
\textit{\ and let }$f$\textit{\ be a composite function given by }%
\begin{equation*}
f=F\left( h_{1}\left( x_{1}\right) \times ...\times h_{n}\left( x_{n}\right)
\right) ,
\end{equation*}%
\textit{where }$h_{1},...,h_{n}$\textit{\ are thrice differentiable and
nonzero functions. Then the Allen matrix }$M\left( f\right) $\textit{\ of }$%
f $\textit{\ is singular if and only if }$f$\textit{\ is one of the
following:}

\textit{(a) }$f=F\left( \gamma e^{\alpha _{1}x_{1}+\alpha _{2}x_{2}}\times
h_{3}\left( x_{3}\right) \times ...\times h_{n}\left( x_{n}\right) \right) ,$%
\textit{\ where }$\gamma ,\alpha _{1},\alpha _{2}$\textit{\ are nonzero
constants;}

\textit{(b) }$f=F\left( \gamma \left( x_{1}+\beta _{1}\right) ^{\alpha
_{1}}\times ...\times \left( x_{n}+\beta _{n}\right) ^{\alpha _{n}}\right) ,$%
\textit{\ where }$\gamma ,\alpha _{i}$\textit{\ are nonzero constants
satisfying }$\alpha _{1}+...+\alpha _{n}=0$\textit{\ and }$\beta _{i}$%
\textit{\ some constants.}

Thus, from Theorem 3.1 and Theorem 4.1, we have the following

\bigskip

\textbf{Corollary 4.2. }\textit{Let }$\left( M^{n},f\right) $\textit{\ be a
homothetical hypersurface in }$%
\mathbb{R}
_{+}^{n+1}$\textit{\ such that at least one of }$\left( \frac{f_{1}^{\prime }%
}{f_{1}}\right) ^{\prime },...,\left( \frac{f_{n}^{\prime }}{f_{n}}\right)
^{\prime }$\textit{\ vanishes. }$\left( M^{n},f\right) $\textit{\ has null
Gauss-Kronocker curvature if and only if the Allen matrix }$H^{B}\left(
f\right) $\textit{\ of }$f$\textit{\ is singular.}

\begin{center}
\textbf{5. A Further Application}
\end{center}

Next result completely classifies the composite functions of the form $%
\left( 4.1\right) $ having constant Hicks elasticity of substitution
property.

\bigskip

\textbf{Theorem 5.1. }\textit{Let }$f\left( \mathbf{x}\right) =F\left(
h_{1}\left( x_{1}\right) \times ...\times h_{n}\left( x_{n}\right) \right) $%
\textit{\ be a twice differentiable production function}$.$\textit{\ Then }$f
$\textit{\ satisfies constant Hicks elasticity of substitution property if
and only if, up to constants, }$f$\textit{\ is one of the following}

\textit{(a) a homothetical generalized Cobb-Douglas production function
given by}%
\begin{equation*}
f=F\left( x^{\alpha _{1}}...x^{\alpha _{n}}\right) ;
\end{equation*}

\textit{(b) a homothetical generalized ACMS production function given by}%
\begin{equation*}
f=F\left( \beta _{1}x_{1}^{\frac{\sigma -1}{\sigma }}+...+\beta _{n}x_{n}^{%
\frac{\sigma -1}{\sigma }}\right) ,\text{ \ \ }\sigma \neq 1;
\end{equation*}

\textit{(c) }$f=F\left( \dprod\limits_{i=1}^{n}\ln \left( x_{i}\right) ^{\mu
_{i}}\right) ,$\textit{\ where }$\mu _{i}$\textit{\ are nonzero constants
for all }$i\in \left\{ 1,...,n\right\} .$

\bigskip

\textbf{Proof. }Let $f$ be a twice differentiable production function given
by%
\begin{equation}
f\left( \mathbf{x}\right) =F\left( h_{1}\left( x_{1}\right) \times ...\times
h_{n}\left( x_{n}\right) \right) .  \tag{5.1}
\end{equation}%
It follows from $\left( 5.1\right) $ that%
\begin{equation}
f_{x_{i}}=\frac{h_{i}^{\prime }}{h_{i}}uF^{\prime },\text{ }f_{x_{i}x_{i}}=%
\frac{h_{i}^{\prime \prime }}{h_{i}}uF^{\prime }+\left( \frac{h_{i}^{\prime }%
}{h_{i}}\right) ^{2}u^{2}F^{\prime \prime }  \tag{5.2}
\end{equation}%
and%
\begin{equation}
f_{x_{i}x_{j}}=\frac{h_{i}^{\prime }h_{j}^{\prime }}{h_{i}h_{j}}u\left(
F^{\prime }+uF^{\prime \prime }\right) ,\text{ }1\leq i\neq j\leq n, 
\tag{5.3}
\end{equation}%
where $u=h_{1}\left( x_{1}\right) \times ...\times h_{n}\left( x_{n}\right)
. $ By substituting $\left( 5.2\right) $ and $\left( 5.3\right) $ into $%
\left( 2.2\right) ,$ we deduce that%
\begin{equation}
\frac{h_{i}^{\prime \prime }h_{i}}{\left( h_{i}^{\prime }\right) ^{2}}+\frac{%
h_{i}}{\sigma x_{i}h_{i}^{\prime }}+\frac{h_{j}^{\prime \prime }h_{j}}{%
\left( h_{j}^{\prime }\right) ^{2}}+\frac{h_{j}}{\sigma x_{j}h_{j}^{\prime }}%
=2.  \tag{5.4}
\end{equation}%
From $\left( 5.4\right) ,$ we have%
\begin{equation}
\frac{h_{i}^{\prime \prime }h_{i}}{\left( h_{i}^{\prime }\right) ^{2}}+\frac{%
h_{i}}{\sigma x_{i}h_{i}^{\prime }}=\zeta _{i},  \tag{5.5}
\end{equation}%
where $\zeta _{i}$ are nonzero constants such that $\zeta _{i}+\zeta _{j}=2$
for $1\leq i\neq j\leq n.$

Now we divide the proof into two separate cases.

\bigskip

\textbf{Case (i): }$\zeta _{i}=1=\zeta _{j}$ for $1\leq i\neq j\leq n.$
After solving $\left( 5.5\right) ,$ we find%
\begin{equation}
h_{i}\left( x_{i}\right) =\left\{ 
\begin{array}{c}
\left. \gamma _{i}x_{i}^{\alpha _{i}}\right. \text{ \ \ \ \ \ \ \ \ \ \ \ \
\ if }\sigma =1 \\ 
\left. e^{\kappa _{i}+\gamma _{i}\left( \frac{\sigma -1}{\sigma }\right)
x_{i}^{\frac{\sigma -1}{\sigma }}}\right. \text{ if }\sigma \neq 1%
\end{array}%
\right.  \tag{5.6}
\end{equation}%
for nonzero constants $\gamma _{i},\alpha _{i}$ and some constant $\kappa
_{i}.$ Combining $\left( 5.1\right) $ and $\left( 5.6\right) $ gives cases
(a) and (b) of the theorem.

\bigskip

\textbf{Case (ii): }$\zeta _{i}\neq 1\neq \zeta _{j}$ for $1\leq i\neq j\leq
n.$ By solving $\left( 5.5\right) ,$ we derive%
\begin{equation}
h_{i}\left( x_{i}\right) =\left\{ 
\begin{array}{c}
\left. \left( \left( 1-\zeta _{i}\right) \kappa _{i}+\ln x_{i}^{\left(
1-\zeta _{i}\right) \eta _{i}}\right) ^{\frac{1}{1-\zeta _{i}}}\right. \text{
\ \ \ \ \ \ \ \ \ \ \ if }\sigma =1 \\ 
\left. \left( \left( \frac{\left( 1-\zeta _{i}\right) \eta _{i}\sigma }{%
\sigma -1}\right) x_{i}^{\frac{\sigma -1}{\sigma }}+\left( 1-\zeta
_{i}\right) \kappa _{i}\right) ^{\frac{1}{1-\zeta _{i}}}\right. \text{ if }%
\sigma \neq 1%
\end{array}%
\right.  \tag{5.7}
\end{equation}%
for nonzero constants $\eta _{i}$ and some constants $\kappa _{i}.$ Now, by
combining $\left( 5.1\right) $ and $\left( 5.7\right) $, up to constants, we
obtain the cases (a) and (c) of the theorem.

Conversely, it is straightforward to verify that each one of cases (a)-(c)
implies that $f$ satisfies constant Hicks elasticity of substitution
property.


\bigskip

\begin{thebibliography}{99}
\bibitem{1} R.G. Allen, J.R. Hicks, A reconsideration of the theory of
value, Pt. II, Economica 1 (1934) 196--219.

\bibitem{2} K. J. Arrow, H. B. Chenery, B. S. Minhas, R. M. Solow,
Capital-labor substitution and economic efficiency, Rev. Econom. Stat. 43
(3) (1961) 225--250.

\bibitem{3} M.E. Aydin, M. Ergut, Homothetic functions with Allen's
perspective and its geometric applications, Kragujevac J. Math. 38 (1)
(2014) 185--194.

\bibitem{4} M.E. Aydin, M. Ergut, Composite functions with Allen
determinants and their applications to production models in economics,
Tamkang J. Math. in press.

\bibitem{5} B. Y. Chen, G. E. V\^{\i}lcu, Geometric classifications of
homogeneous production functions, Appl. Math. Comput. 225 (2013) 345--351.

\bibitem{6} B. Y. Chen, S. Decu and L. Verstraelen, Notes on isotropic
geometry of production models, Kragujevac J. Math. 38 (1) (2014) 23--33.

\bibitem{7} B.-Y. Chen, Geometry of submanifolds, M. Dekker, New York, 1973.

\bibitem{8} B.-Y. Chen, On some geometric properties of h-homogeneous
production function in microeconomics,\ Kragujevac J. Math. 35 (3) (2011)
343--357.

\bibitem{9} B. Y. Chen, On some geometric properties of quasi-sum production
models, J. Math. Anal. Appl. 392 (2012) 192--199.

\bibitem{10} B.-Y. Chen, An explicit formula of Hessian determinants of
composite functions and its applications, Kragujevac J. Math. 36 (1) (2012)
27--39.

\bibitem{11} B.-Y. Chen, Geometry of quasi-sum production functions with
constant elasticity of substitution property, J. Adv. Math. Stud. 5 (2)
(2012) 90--97.

\bibitem{12} B.-Y. Chen, Classification of homothetic functions with
constant elasticity of substitution and its geometric applications, Int.
Elect. J. Geo. 5 (2)\ (2012) 67--78.

\bibitem{13} B.-Y. Chen, Classification of h-homogeneous production
functions with constant elasticity of substitution, Tamkang J. Math. 43
(2012) 321--328.

\bibitem{14} B. Y. Chen,\ A note on homogeneous production models,
Kragujevac J. Math. 36 (1) (2012), 41--43.

\bibitem{15} B.-Y. Chen,\ Solutions to homogeneous Monge-Ampere equations of
homothetic functions and their applications to production models in
ecenomics, J. Math. Anal. Appl. 411 (2014) 223--229.

\bibitem{16} C. W. Cobb, P. H. Douglas, A theory of production, Amer.
Econom. Rev. 18 (1928) 139--165.

\bibitem{17} M.P. do Carmo, Differential Geometry of Curves and Surfaces,
Prentice Hall, Englewood Cliffs, NJ, 1976.

\bibitem{18} W. Goemans, I. Van de Woestyne, Translation and homothetical
lightlike hypersurfaces of semi-Euclidean space, Kuwait J. Sci. Eng. 38 (2A)
(2011) 35-42.

\bibitem{19} A. Gray, Modern differential geometry of curves and surfaces
with mathematica,\ CRC Press LLC, 1998.

\bibitem{20} L. Jiu, H. Sun, On minimal homothetical hypersurfaces, Colloq.
Math. 109 (2007) 239--249.

\bibitem{21} P. O. Linderberg, E. A. Eriksson and L. G. Mattsson, Homothetic
functions revisited, Econ. Theory 19 (2002) 417-427.

\bibitem{22} L. Losonczi, Production functions having the CES property, Acta
Math. Acad. Paedagog. Nyh\'{a}i. (N.S.) 26 (1) (2010) 113--125.

\bibitem{23} S.K. Mishra, A brief history of production functions, IUP J.
Manage. Econom.\ 8 (4) (2010) 6--34.

\bibitem{24} D. Saglam, A. Sabuncuoglu, Geodesic curves of the minimal
homothetical hypersurfaces in the semi-Euclidean space, Diff. Geom. Dyn.
Syst. 10\ (2008) 275-287.

\bibitem{25} D. Saglam, A. Sabuncuoglu, Minimal homothetical lightlike
(degenerate) hypersurfaces of semi-Euclidean spaces, Kuwait J. Sci. Eng. 38
(1A) (2011), 1-14.

\bibitem{26} H. Uzawa,\ Production functions with constant elasticities of
substitution, The Review of Economic Studies\ 29 (4) (1962) 291-299.

\bibitem{27} I. Van de Woestyne, Minimal homothetical hypersurfaces of a
semi-Euclidian space, Results in Mathematics 27\ (1995) 333--342.

\bibitem{28} G.E. V\^{\i}lcu,\ A geometric perspective on the generalized
Cobb--Douglas production functions, Appl. Math. Lett. 24\ (2011) 777--783.

\bibitem{29} A. D. Vilcu, G. E. Vilcu, On some geometric properties of the
generalized CES production functions, Appl. Math. Comput. 218 (2011)
124--129.

\bibitem{30} X. Wang, Y. Fu, Some characterizations of the Cobb-Douglas and
CES production functions in microeconomics, Abstract Appl. Anal.
doi.org/10.1155/2013/761832.
\end{thebibliography}
\end{document}